# STRUCTURE THEOREMS FOR EMBEDDED DISKS WITH MEAN CURVATURE BOUNDED IN $L^P$

GIUSEPPE TINAGLIA

ABSTRACT. After appropriate normalizations an embedded disk whose second fundamental form has large norm contains a multi-valued graph, provided the $L^p$ norm of the mean curvature is sufficiently small. This generalizes to non-minimal surfaces a well known result of Colding and Minicozzi.

## Introduction

In [7] Colding and Minicozzi proved that a minimal disk embedded in $\mathbb{R}^3$ whose Gaussian curvature is large at a point contains a multi-valued graph around that point. This means that, locally, the disk looks like a piece of a suitably scaled helicoid (see Figure 1). This was later generalized in [20] to the constant mean curvature case. The structure theorem in [7] has been used as a key ingredient in the Colding-Minicozzi series of papers [3, 6, 7, 8, 9] dealing with the geometry of embedded minimal surfaces of fixed genus. Moreover, the new ideas provided by their recent work have been applied to solve several longstanding problems in the field of minimal surfaces; see for instance [2, 4, 13, 14].

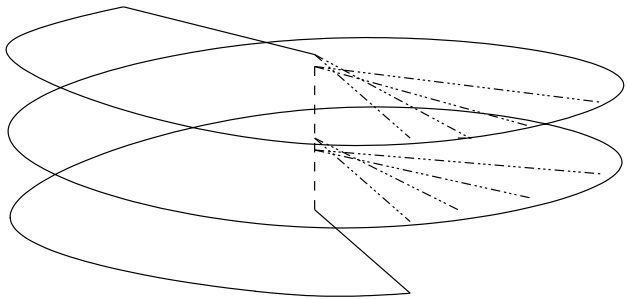

FIGURE 1. Half of the helicoid

In this paper we discuss the geometry of disks embedded in $\mathbb{R}^3$ for which the $L^p$ norm of the mean curvature, $\|H\|_{L^p}$, is suitably bounded. We point out that $\|H\|_{L^p}$ is a natural quantity to consider. It appears, for instance, in such classical results as the monotonicity formulae [19] and related applications [10, 12]. Loosely speaking, the principle established in our main result is that an embedded disk whose second fundamental form is large at the origin, and whose $L^q$-norm of the gradient of the mean curvature is bounded, $q > 2$, must contain a multi-valued graph if the $L^p$-norm of the mean curvature is suitably small. Below is a simplified version of the main theorem.

**Theorem 0.1.** *Given $N \in \mathbb{Z}_+$, $T \geq 0$, $q > 2$ and $p \geq 1$ there exist $C_1 = C_1(N) > 0$, $C_2 = C_2(N, T, p, q) > 0$, and $l = l(N, p) > 1$ such that the following holds.
If $\Sigma \subset \mathbb{R}^3$ is an embedded disk with $0 \in \Sigma \subset B_l(0)$, $\partial \Sigma \subset \partial B_l(0)$, $\|H\|_{L^p} \leq C_2$, $\|\nabla H\|_{L^q} \leq T$, and*
$$\sup_{\Sigma \cap B_l(0)} |A_\Sigma| \leq 2C_1 = 2|A_\Sigma|(0),$$
*then $\Sigma \cap B_1(0)$ contains an $N$-valued graph that forms around the origin.*

Here $B_l(0)$ is the euclidean ball of radius $l$ centered at the origin. We recall that if $\Sigma$ is a surface and $k_1$ and $k_2$ are its principal curvatures, then the mean curvature is $H = \frac{k_1+k_2}{2}$. The norm of the second fundamental form is $|A_\Sigma| = \sqrt{k_1^2 + k_2^2}$ and the Gaussian curvature is $K_\Sigma = k_1 k_2$. A precise definition of an $N$-valued graph as well as a finer quantitative version of Theorem 0.1 is to be found in Section 4.

Our new generalization of the Colding-Minicozzi structure theorem is intended as a first step towards classifying the singularities of the limit of a sequence of embedded disks with $\|H\|_{L^p}$ bounded. While this problem has been successfully studied for minimal disks [3, 5, 9, 16], it remains unsolved in this more general setting. In fact, if the norm of the second fundamental form of these disks is uniformly bounded, they do converge to a well defined surface, although not necessarily an embedded one. The main objective is therefore to understand what happens as the norm of the second fundamental form is arbitrarily large. As is the case of minimal surfaces, the answer to this question will surely provide new tools for the study of the global properties of surfaces with $\|H\|_{L^p}$ bounded.

The proof of Theorem 0.1 uses a new compactness argument which, somewhat unexpectedly, does not require a bound on the area.

## Acknowledgments

It is a pleasure to thank the geometry group at Stanford University, where this work was initiated while the author was a post-doctoral fellow.

## 1. Minimal Surfaces

1. **Definition.** Let $\Sigma \subset \mathbb{R}^3$ be a 2-dimensional smooth orientable surface (possibly with boundary) with unit normal $N_\Sigma$. Given a function $\phi$ in the space $C_0^\infty(\Sigma)$ of infinitely differentiable, compactly supported functions on $\Sigma$, consider the one-parameter variation
$$\Sigma_{t,\phi} = \{x + t\phi(x)N_\Sigma(x) \mid x \in \Sigma\}$$
and let $A(t)$ be the area functional,
$$A(t) = Area(\Sigma_{t,\phi}).$$
The so-called first variation formula is the equation

$$\text{(1.1)} \qquad \frac{dA}{dt}(0) = -2 \int_\Sigma \phi H,$$

where $H$ is the mean curvature of $\Sigma$. If a surface $\Sigma$ is given as graph of a function $u$ then

$$2H = \text{div}\left(\frac{\nabla u}{\sqrt{1+|\nabla u|^2}}\right).$$

A surface is said to be a *minimal* surface if it is a critical point for the area functional (see [17]); by equation (1.1) this is equivalent to require that the mean curvature be identically zero. Examples of minimal surfaces are planes, the helicoid and the catenoid.

2. **Limits of minimal surfaces.** In this section we discuss limits of minimal surfaces. Some of this material is covered in great detail (including proofs) in [18, Section 4].

Let $\Sigma$ be a surface in $\mathbb{R}^3$ and let $T_p\Sigma$ denote its tangent plane at $p$. Given $p \in \Sigma$ and $r > 0$ we label by

$$D(p, r) = \{p + v | v \in T_\Sigma, |v| < r\}$$

the tangent disk of radius $r$. We use $W(p, r)$ to stand for the infinite solid cylinder of radius $r$ around the affine normal line at $p$ is noted by,

$$W(p, r) = \{q + tN(p) | q \in D(p, r), t \in \mathbb{R}\}.$$

Inside $W(p, r)$ and for $\varepsilon > 0$, we have the compact slice

$$W(p, r, \varepsilon) = \{q + tN(p) | q \in D(p, r), |t| < \varepsilon\}.$$

**Definition 1.2.** *Let $\Sigma_n$ be a sequence of surfaces properly embedded in an open set $O \subset \mathbb{R}^3$. We say that $\Sigma_n$ converges $C^k$ with finite multiplicity to a surface $\Sigma_\infty$ on compact sets if for any $K \subset O$ compact there exist $r, \varepsilon > 0$ such that*
  (1) *For any $p \in \Sigma_\infty \cap K$, $\Sigma_\infty \cap W(p, r, \varepsilon)$ can be represented as the graph of a function $u \colon D(p, r) \to \mathbb{R}$.*
  (2) *For sufficiently large $n$, $\Sigma_n \cap K \cap W(p, r, \varepsilon)$ consists of a finite number of graphs (independent of $n$) over $D(p, r)$ which converge to $u$ in the usual $C^k$ topology.*

Given a sequence of subsets $\{F_n\}_n$ in an open domain $O$, its *accumulation set* is defined by $\{p \in O | \exists p_n \in F_n \text{ with } p_n \longrightarrow p\}$

The next two theorems are standard compactness theorems for sequences of properly embedded minimal surfaces. The first assumes a uniform bound on the area and on the norm of the second fundamental form, while the second one is slightly more general and does not assume a bound on the area. Their proofs are similar in nature. We will only sketch the proof of Theorem 1.4, pointing out where to use the uniform bound on the area to prove Theorem 1.3.

**Theorem 1.3.** *Let $\Sigma_n$ be a sequence of minimal surfaces properly embedded in an open set $O \subset \mathbb{R}^3$. Suppose that $\Sigma_n$ has an accumulation point in $O$ and that there exist $C_1$ and $C_2$ such that*
$$Area(\Sigma_n) < C_1, \quad \sup_{\Sigma_n} |A_n| < C_2$$
*uniformly in $n$. Then, there exists a subsequence $\Sigma_{n_k}$ and a minimal surface $\Sigma$ properly embedded in $O$ such that $\Sigma_{n_k}$ converges smoothly with finite multiplicity to $\Sigma$ on compact subsets of $O$.*

**Theorem 1.4.** *Let $\Sigma_n$ be a sequence of minimal surfaces properly embedded in an open set $O \subset \mathbb{R}^3$. Suppose that there exists a sequence $p_n \in \Sigma_n$ converging to a point $p \in O$ and that*
$$\sup_{\Sigma_n} |A_n| < C$$
*uniformly in $N$. Then, there exists a subsequence $\Sigma_{n_k}$ and a connected minimal surface $\Sigma$ in $O$ satisfying*
  (1) *$\Sigma$ is contained in the accumulation set of $\{\Sigma_n\}_n$;*
  (2) *$p \in \Sigma$ and $|A|(p) = \lim |A_n|(p_n)$;*
  (3) *$\Sigma$ is embedded in $O$;*
  (4) *Any divergent path in $\Sigma$ either diverges in $O$ or has infinite length.*

*Sketch of the proof of 1, 2 and 3.* As $p_n$ accumulates at $p \in O$, the uniform bound on the second fundamental form implies that there exists $r > 0$ such that for $\varepsilon$ small, $W(p_n, r, \varepsilon) \cap \Sigma_n$ consists of a collection of graphs, given by the functions $u_n^k \colon D(p_n, r) \to \mathbb{R}$ (regarding Theorem 1.3 a uniform bound on the area would give a bound for the number of graphs which is independent of $n$). After going to a subsequence we can assume that $T_{p_n}\Sigma_n$ converges to a plane $\pi$ and that the functions $u_n^k$ are defined over $\pi$. Moreover, $|u_n^k|$, $|\nabla u_n^k|$ and $|\nabla^2 u_n^k|$ are uniformly bounded. Since $u_n$ is a minimal graph, thanks to the minimal graph equation we have uniform $C^{2,\alpha}$ estimates for $u_n^k$. In this situation, Arzela-Ascoli's Theorem implies that a subsequence of $u_n^k$ converges $C^2$ to a function $u$. Due to the $C^2$ convergence, $u$ is also a minimal graph. An analytic continuation argument allows us to construct a subsequence $\Sigma_k$ and a maximal sheet $\Sigma$ in the accumulation set of $\Sigma_k$ which extends the graph $u$. By construction, $\Sigma$ satisfies items (1) and (2). That (3) also holds can be seen as follows. The surface $\Sigma$ must be embedded because transverse self-intersections of it would give rise to transverse self intersection of $\Sigma_n$ for $n$ large and tangential self-intersections would contradict the maximum principle. □

## 2. Compactness Theorems

In this section we prove more general compactness theorems which will be used in the proof of the structure theorem. We start by extending Theorem 1.3 and 1.4 to more general surfaces. Before doing that, we need to establish some notation. Let $f$ be a function defined over $\Sigma$ and let $p \geq 1$, then $\|f\|_{L^p(\Sigma)}$ is the $L^p$ norm of $f$, while $\|f\|_\alpha$ is the $C^\alpha$ norm.

**Theorem 2.5.** Let $\Sigma_n$ be a sequence of smooth surfaces properly embedded in an open set $O \subset \mathbb{R}^3$. Suppose that $\Sigma_n$ has an accumulation point in $O$ and that there exist $C_1$, $C_2$, and $T$ such that
$$Area(\Sigma_n) < C_1, \sup_{\Sigma_n} |A_n| < C_2 \text{ and } \|H_n\|_\alpha < T \text{ uniformly in } n.$$
Suppose also that $\|H_n\|_{L^p(\Sigma_n)}$ converges to zero for some $p \geq 1$. Then, there exists a subsequence $\Sigma_{n_k}$ and a minimal surface $\Sigma$ properly embedded in $O$ such that $\Sigma_{n_k}$ converges $C^2$ with finite multiplicity to $\Sigma$ on compact sets of $O$.

**Theorem 2.6.** Let $\Sigma_n$ be a sequence of surfaces properly embedded in an open set $O \subset \mathbb{R}^3$. Suppose that there exists a sequence $p_n \in \Sigma_n$ converging to a point $p \in O$ and that there exist $C$ and $T$ such that
$$\sup_{\Sigma_n} |A_n| < C \text{ and } \|H_n\|_\alpha < T \text{ uniformly in } n.$$
Suppose also that $\|H_n\|_{L^p(\Sigma_n)}$ converges to zero for some $p \geq 1$. Then, there exists a subsequence $\Sigma_{n_k}$ and a connected minimal surface $\Sigma$ in $O$ satisfying

(1) $\Sigma$ is contained in the accumulation set of $\Sigma_n$.
(2) $p \in \Sigma$ and $|A|(p) = \lim |A_n|(p_n)$.
(3) $\Sigma$ is embedded in $O$.
(4) Any divergent path in $\Sigma$ either diverges in $O$ or has infinite length.

Their proofs are a slight modification of the proofs of Theorem 1.3 and 1.4. As before, the uniform bound on the second fundamental form implies that $W(p_n, r, \varepsilon) \cap \Sigma_n$ consists of a collection of graphs. In the proofs of Theorem 1.3 and 1.4 we needed the surfaces to be minimal in order to obtain uniform $C^{2,\alpha}$ estimates for these graphs. In fact, it is known that in order to obtain $C^{2,\alpha}$ estimates it suffices that $\|H_n\|_\alpha$ is uniformly bounded. Once uniform $C^{2,\alpha}$ estimates are obtained we can apply the Arzela-Ascoli theorem to extract a subsequence of graphs which converges $C^2$ to a graph. The fact that $\|H_n\|_{L^p(\Sigma_n)}$ converges to zero is ultimately used to show that the limit graph is minimal.

We remark that using the Sobolev embedding theorem, assuming a bound on $\|\nabla H_n\|_{L^q(\Sigma_n)}$, $q > 2$, implies a bound on $\|H\|_\alpha$. We could therefore restate Theorem 1.3 and 1.4 replacing the uniform bound on $\|H\|_\alpha$ with a uniform bound on $\|\nabla H_n\|_{L^q(\Sigma_n)}$, $q > 2$. Furthermore, if $\|H\|_\alpha$ is not bounded uniformly, although we are not able to extract a subsequence of graphs which converges $C^2$ to a minimal graph, it is still possible to extract a subsequence of graphs which converges $C^1$ to a minimal graph. Moreover, an upper bound on the norm of the second fundamental form of the limit is still valid. In other words, the following theorem follows.

**Theorem 2.7.** Let $\Sigma_n$ be a sequence of surfaces properly embedded in an open set $O \subset \mathbb{R}^3$. Suppose that there exists a sequence $p_n \in \Sigma_n$ converging to a point $p \in O$ and that there exists $C$ such that
$$\sup_{\Sigma_n} |A_n| < C, \text{ uniformly in } n.$$

Suppose also that $\|H_n\|_{L^p(\Sigma_n)}$ converges to zero for some $p \geq 1$. Then, there exists a subsequence $\Sigma_{n_k}$ and a connected minimal surface $\Sigma$ in $O$ satisfying
  (1) $\Sigma$ is contained in the accumulation set of $\Sigma_n$.
  (2) $\sup_\Sigma |A| < C$.
  (3) $\Sigma$ is embedded in $O$.
  (4) Any divergent path in $\Sigma$ either diverges in $O$ or has infinite length.

In the next theorem we use Theorem 2.7 to describe more accurately the accumulation set of a sequence of surfaces whose $\|H_n\|_{L^p}$ converges to zero.

**Theorem 2.8.** *Let $\Sigma_n$ be a sequence of compact surfaces embedded in $B_n(0)$ such that $0 \in \Sigma \subset B_n(0)$, $\partial\Sigma \subset \partial B_n(0)$. Suppose that there exists a constant $C$ such that $\sup_{\Sigma_n} |A_n| < C$ and that $\|H_n\|_{L^p}$ converges to zero for some $p \geq 1$ as $n$ goes to infinity. Then, up to a subsequence, the accumulation set of $\Sigma_n$ is non-empty and it consists either of a connected complete properly embedded minimal surface $\Sigma$ or of a collection of parallel planes.*

*Moreover, $\Sigma_n^1$ converges $C^1$ with multiplicity one to $\Sigma^1$, where $\Sigma_n^1$ is the connected component of $\Sigma_n \cap B_1(0)$ which contains the origin and $\Sigma^1$ is, depending on the accumulation set, the connected component of $\Sigma \cap B_1(0)$ which contains the origin or a unit disk centered at the origin.*

*Proof.* Supposing the conclusion fails, let $\Sigma_n$ be a sequence of compact surfaces embedded in $B_n(0)$ such that $0 \in \Sigma_n \subset B_n(0)$, $\partial\Sigma_n \subset \partial B_n(0)$ and $\|H_n\|_{L^p} < \frac{1}{n}$. Theorem 2.7 implies that there exists a complete connected embedded minimal surface, $\Sigma$, which contains the origin and is contained in the accumulation set of $\Sigma_n$. Furthermore, $\Sigma$ has bounded second fundamental form and therefore it is properly embedded (see [15]). If $\Sigma_n$ has another accumulation point which is not in $\Sigma$, then the same argument shows that there exists another complete connected properly embedded minimal surface, $\Sigma'$, which is contained in the accumulation set of $\Sigma_n$ and it is disjoint from $\Sigma$. The results in [1, 11, 21] imply that they must be parallel planes.

Let $\varepsilon > 0$ and let $TN(\varepsilon)$ be an embedded tubular neighborhood of $\Sigma^1$ of size $\varepsilon$. Choose $r$ and $\varepsilon$, $r > 2\varepsilon > 0$ such that for any $p \in \Sigma$, $W(p,r,\varepsilon) \cap \Sigma_n^1$ consists of a collection of graphs. Let $u_0$ be the minimal graph over $D(0,r)$ which locally represents $\Sigma^1$, and let $u_0^n$ be the graph in $W(p,r,\varepsilon) \cap \Sigma_1^n$ containing the origin. From the way $\Sigma$ has been obtained, $u_0^n$ converges $C^1$ to $u_0$. For any $q \in \partial W(0,r,\varepsilon)$ if we let $u_q$ represents $W(q,r,\varepsilon) \cap \Sigma^1$ and $u_q^n$ represents the connected component of $W(q,r,\varepsilon) \cap \Sigma_n^1$ which intersects $u_0^n$, then we can assume that $u_q^n$ converges $C^1$ to $u_q$. Since $\Sigma$ is properly embedded, $\Sigma^1$ is compact. After finitely many steps it is possible to analytically continue $u_0^n$ to get a one-sheeted cover of $\Sigma_1$. □

Loosely speaking, in the next theorems we describe the geometry away form the boundary of a surface whose $L^p$ norm of the mean curvature is small.

**Theorem 2.9.** *Given $C > 0$, $p \geq 1$ there exist $R = R(C,p) > 2$, $\varepsilon = \varepsilon(C,p) > 0$ such that the following holds. Let $\Sigma$ be a compact surface embedded in $B_R(0)$ such* that $0 \in \Sigma \subset B_R(0)$, $\partial \Sigma \subset \partial B_R(0)$, $\sup_\Sigma |A| < C$ and $\|H\|_{L^p} < \frac{1}{R}$ then $\Sigma^1$ has an embedded tubular neighborhood of size $\varepsilon$.

*Proof.* If not, let $\Sigma_n$ be a sequence of compact surfaces embedded in $B_n(0)$ such that $0 \in \Sigma_n \subset B_n(0)$, $\partial \Sigma \subset \partial B_n(0)$, $\|H_n\|_{L^p} < \frac{1}{n}$ and such that the size of the largest embedded tubular neighborhood of $\Sigma_n^1$ is going to zero. Theorem 2.8 implies that, after going to a subsequence, $\Sigma_n^1$ converges $C^1$ with multiplicity one to a compact properly embedded minimal surface. However, the fact that the size of the tubular neighborhood is going to zero clearly contradicts the multiplicity one convergence. $\square$

An easy consequence of Theorem 2.9 is the existence of some upper bound for the area of a compact embedded surface depending solely on the upper bounds for the $L^p$ norm of the mean curvature and second fundamental form. Although the area of such a surface is not necessarily uniformly bounded, we show that there exists a uniform upper bound for the area of connected pieces which are sufficiently away from the boundary.

**Corollary 2.10.** *Given $C > 0$, $p \geq 1$ there exist $K = K(C, p) > 0$ and $R = R(C, p) > 2$ such that the following holds. Let $0 \in \Sigma$ be a compact embedded surface such that $\Sigma \subset B_R(0)$, $\partial \Sigma \subset \partial B_R(0)$, $\|H\|_{L^p} < \frac{1}{R}$ and $\sup_\Sigma |A| < C$ then the area of $\Sigma^1$ is bounded by $K$.*

In the following compactness theorem we prove that if the elements of the sequence in Theorem 2.8 are embedded disks, so is the limit.

**Theorem 2.11.** *Given $C > 0$, $p \geq 1$ there exists $R = R(C, p) > 2$ such that the following holds. Let $\Sigma_n$ be a surface embedded in $B_R(0)$ such that $0 \in \Sigma \subset B_R(0)$, $\partial \Sigma \subset \partial B_R(0)$ and $\sup_{\Sigma_n} |A_n| < C$. Suppose also that $\|H\|_{L^p}$ converges to zero as $n$ goes to infinity and that $\Sigma_n$ is simply connected. Then, up to a subsequence, $\Sigma_n^1$ converges $C^1$ to a properly embedded minimal disk $\Sigma^1$.*

*Proof.* In light of Theorem 2.8, all that needs to be showed is that $\Sigma_n^1$ is a disk. This will be discussed in the next section. $\square$

## 3. Weak Convex Hull Properties

In this section we prove a weak convex hull property for surfaces with bounded $L^p$ norm of the mean curvature. Loosely speaking, we say that a surface satisfies a weak convex hull property if when the boundary is contained in a ball $B$, then the entire surface has to be contained in $B$.

In the next lemma we prove that if $\Sigma$ is a surface contained in a compact set, its second fundamental form is bounded and its boundary is contained in a certain ball then, if the $L^p$ norm of the mean curvature is small enough, the surface cannot stretch too far outside the ball. The proof is by contradiction and uses a compactness argument. The idea is that after taking a convergent subsequence, since the limit minimal surface satisfies a convex hull property, an analogous property has to be satisfied by the elements in the sequence. Notice that we need the elements of the sequence to be contained in a compact set otherwise one could take a sequence of spheres with radii going to infinity. The $L^p$ norm of the mean curvature of these spheres converges to zero, $p > 2$, but they do not satisfy any weak convex hull property.

**Lemma 3.12.** *Given $l > 1$, $\varepsilon > 0$ and $1 \leq p < \infty$ there exists an $n = n(l, \varepsilon, p) > 0$ such that the following holds. Suppose $\Sigma$ is a compact surface such that $\Sigma \subset B_l(0)$, $\partial \Sigma \subset \partial B_l(0)$, $\|H\|_{L^p} < \frac{1}{n}$ and $\sup_\Sigma |A| < C$ and let $\Sigma^1 \subset \Sigma$ be a surface such that $\partial \Sigma^1 \subset B_1(0)$. Then $\Sigma^1 \subset B_{1+\varepsilon}(0)$.*

*Proof.* Assume, by way of contradiction, that there exists a sequence of $\Sigma_n$ and $\Sigma_n^1 \subset \Sigma_n$ such that $\Sigma_n \subset B_l(0)$, $\|H_n\|_{L^p} < \frac{1}{n}$, $\sup_{\Sigma_n} |A_n| < C$, $\partial \Sigma_n^1 \subset B_1(0)$ and $\Sigma_n^1 \not\subset B_{1+\varepsilon}(0)$. Let $p_n \in \Sigma_n^1$ such that

$$(3.2) \qquad l \geq |p_n| = \max_{q \in \Sigma_n^1 \cap \mathbb{R}^3 \setminus B_1(0)} |q| > 1 + \varepsilon.$$

After going to a subsequence we can assume that $p_n$ converges to a point $p \in B_l(0) \setminus B_{1+\varepsilon}(0)$. Consider $\delta \leq \frac{\varepsilon}{2}$ such that the connected component of $W(p_n, \delta, \delta)$ that contains $p_n$ consists of a graph over $D(p_n, \delta)$. In particular, after going to a subsequence, the graph containing $p_n$ would converge $C^1$ to a minimal graph which is tangent to $B_{|p|}(0)$ and contained inside its convex side. This contradicts the maximum principle and proves the theorem. □

In the case when $p = \infty$ the weak convex hull property does not require a bound on the second fundamental form and can be proved without using a compactness argument.

**Lemma 3.13.** *Fix $l > 1$ and let $\Sigma$ be a compact embedded surface such that $\Sigma \subset B_l(0)$, $\partial \Sigma \subset \partial B_l(0)$, $\sup_\Sigma |H| < \frac{1}{2l}$. Let $\Sigma' \subset \Sigma$ be a surface such that $\partial \Sigma' \subset B_r(p)$, $r > 0$, then $\Sigma' \subset B_r(p)$.*

*Proof.* If $\Sigma'$ is not contained in $B_r(p)$ then there exists an $R$, $r < R < 2l$, such that $\Sigma'$ is contained inside $B_R(p)$ and it is tangent to its boundary. Let $k_1(q)$ and $k_2(q)$ be the principal curvatures at $q$. Clearly, $k_1(q)$ and $k_2(q)$ have the same sign, and $|k_i(q)| \geq \frac{1}{R} > \frac{1}{2l}$. Consequently, $|H(q)| \geq \frac{1}{2l}$. This contradicts the assumption and proves the lemma. □

## 4. STRUCTURE THEOREM

In this section we prove the structure theorem for embedded disks with bounded $L^p$ norm of the mean curvature, $p \geq 1$, and bounded $L^q$ norm of the gradient of the mean curvature, $q > 2$. For simplicity we are going to state the theorems when $p = \infty$ and assuming a bound on $[H]_\alpha := \sup_\Sigma \frac{|H(x) - H(y)|}{dist_\Sigma(x,y)}$.

**Definition 4.14** (Multi-valued graph (see [7])). *Let $D_r$ be the disk in the plane centered at the origin and of radius $r$ and let $\mathcal{P}$ be the universal cover of the punctured plane $\mathbb{C} \setminus 0$ with global coordinates $(\rho, \theta)$ so $\rho > 0$ and $\theta \in \mathbb{R}$. An $N$-valued graph* of a function $u$ on the annulus $D_s\setminus D_r$ is a single valued graph over $\{(\rho,\theta)|r \leq \rho \leq s, |\theta| \leq N\pi\}$.

The prototypical example of a surface which contains a multi-valued graph is the helicoid, Fig. 1. A parametrization of the helicoid that illustrates the existence of such an $N$-valued graph is the following

$$(s\sin t, s\cos t, t) \qquad \text{where } (s,t) \in \mathbb{R}^2.$$

It is easy to see that it contains the $N$-valued graph $\phi$ defined by

$$\phi(\rho,\theta) = \theta \qquad \text{where } (\rho,\theta) \in \mathbb{R}^+\setminus 0 \times [-N\pi, N\pi].$$

In their fundamental work Colding and Minicozzi proved the following structure theorem for minimal disks.

**Theorem 4.15.** [7, **Theorem 0.4.**] *Given $N \in \mathbb{Z}_+$, $\omega > 1$ and $\varepsilon > 0$, there exist $C = C(N,\omega,\varepsilon) > 0$ such that the following holds.*

*Let $0 \in \Sigma \subset B_R \subset \mathbb{R}^3$ be a compact embedded minimal disk such that $\partial\Sigma \subset \partial B_R$. If*

$$\sup_{\Sigma \cap B_{r_0}} |A|^2 \leq 4C^2 r_0^{-2} \text{ and } |A|^2(0) = C^2 r_0^{-2}$$

*for some $0 < r_0 < R$, then there exists $\overline{R} < \frac{r_0}{\omega}$ and (after a rotation) an $N$-valued graph $\Sigma_g \subset \Sigma$ over $D_{\omega\overline{R}}\setminus D_{\overline{R}}$ with gradient $\leq \varepsilon$ and $\text{dist}_\Sigma(0,\Sigma_g) \leq 4\overline{R}$.*

The result below extends Theorem 4.15 to the non-minimal setting.

**Theorem 4.16.** *Given $N \in \mathbb{Z}_+$, $\omega > 1$, $\varepsilon > 0$, $p \geq 1$, and $T > 0$ there exist $C_1 = C_1(N,\omega,\varepsilon) > 0$, $C_2 = C_2(N,\omega,\varepsilon,T,p) > 0$, and $l = l(N,\omega,\varepsilon,p) > 1$ such that the following holds. If $\Sigma \subset \mathbb{R}^3$ is a compact embedded disk with $0 \in \Sigma \subset B_{r_0 l}(0)$, $\partial\Sigma \subset \partial B_{r_0 l}(0)$,*

$$\sup_{\Sigma \cap B_{r_0 l}(0)} |A|^2 \leq 4C_1^2 r_0^{-2} \text{ and } |A|^2(0) = C_1^2 r_0^{-2},$$

(4.3) $$r_0|H| \leq \min(C_2, \frac{1}{2l}) \text{ and } r_0^{1+\alpha}[H]_\alpha \leq T$$

*for some $r_0 > 0$, then there exists $\overline{R} < \frac{r_0}{\omega}$ and (after a rotation) an $N$-valued graph $\Sigma_g \subset \Sigma$ over $D_{\omega\overline{R}}\setminus D_{\overline{R}}$ with gradient $\leq \varepsilon$ and $\text{dist}_\Sigma(0,\Sigma_g) \leq 4\overline{R}$.*

*Proof.* Theorem 4.16 will follow by rescaling after we prove it for $r_0 = 1$. Assuming $r_0 = 1$ the hypotheses become $0 \in \Sigma \subset B_l(0)$, $\partial\Sigma \subset \partial B_l(0)$,

$$\sup_{\Sigma \cap B_{\bar{l}}(0)} |A|^2 \leq 4C_1^2 \text{ and } |A|^2(0) = C_1^2,$$

$$|H| \leq \min(C_2, \frac{1}{2l}) \text{ and } [H]_\alpha \leq T.$$

We have to prove that fixed $T > 0$ there exists $C_2$ such that if the above is true, then $\Sigma$ contains a multi-valued graph in the ball of radius 1. The proof is by contradiction and uses a compactness argument.

Assuming that the theorem is false, let $C_1$ be as big as given by Theorem 4.15, $l$ as given by Theorem 2.11, and $\Sigma_n$ a sequence of embedded disks satisfying the hypotheses of the statement that does not contain a multi-valued graph and with $|H|$ less than $\frac{1}{n}$. As $n$ goes to infinity, Theorem 2.11 gives that, up to a subsequence, $\Sigma_n^1$ (the connected component of $\Sigma_n \cap B_1(0)$ containing the origin) converges $C^2$ with multiplicity one to a minimal disk $\Sigma^1$. The minimal disk containing the origin satisfies the hypotheses of Theorem 4.15 and therefore it contains a multi-valued graph. Since the limit contains a multi-valued graph, $\Sigma_n^1$ must also contain a multi-valued graph for $n$ large. This gives a contradiction and proves the theorem. $\square$

Notice that the $C^2$ convergence guarantees that not only does $\Sigma_n$ contain an $N$-valued graph, but that the properties of this graph, such as the upper bound on the gradient, are preserved.

When the mean curvature is bounded in $L^\infty$ we can prove the next two corollaries. For simplicity we will not state them in full generality as we did for Theorem 4.16. The general versions can be easily obtained using a rescaling argument.

In the next corollary we prove that if the second fundamental form of a compact embedded disk at the origin is bigger than what it is necessary to prove existence of an $N$-valued graph, it is almost its maximum, and the disk satisfies (4.3) in Theorem 4.16, then the disk must contain a multi-valued graph, possibly on a smaller scale.

**Corollary 4.17.** *Given $N \in \mathbb{Z}_+$, $\omega > 1$, $\varepsilon > 0$, and $T > 0$ there exist $C_1 = C_1(N, \omega, \varepsilon) > 0$, $C_2 = C_2(N, \omega, \varepsilon) > 0$, and $l_1 = l_1(N, \omega, \varepsilon) > 1$ such that the following holds.*

*If $\Sigma \subset \mathbb{R}^3$ is a compact embedded disk with $0 \in \Sigma \subset B_l(0)$, $\partial \Sigma \subset \partial B_l(0)$,*

$$\sup_{\Sigma \cap B_l(0)} |A|^2 \leq 4(C + \beta)^2 \text{ and } |A|^2(0) = (C + \beta)^2,$$

$$|H| \leq \min(C_2, \frac{1}{2l}) \text{ and } [H]_\alpha \leq T$$

*for some $\alpha > 0$, $l > l_1$ then there exists $\overline{R} < \frac{C}{\omega(C+\beta)}$ and (after a rotation) an $N$-valued graph $\Sigma_g \subset \Sigma$ over $D_{\omega \overline{R}} \setminus D_{\overline{R}}$ with gradient $\leq \varepsilon$.*

*Proof.* Consider the rescaled surface $\Sigma' = \frac{C+\beta}{C} \Sigma$ and let $\Sigma''$ be the connected component of $\Sigma' \cap B_l(0)$ that contains the origin. Thanks to the weak convex hull property $\Sigma''$ is still a disk, and since $\frac{C+\beta}{C} > 1$, it satisfies the hypothesis of Theorem 4.16. It follows that there exists $\overline{R} < \frac{1}{\omega}$ and (after a rotation) an $N$-valued graph $\Sigma_g \subset \Sigma''$ over $D_{\omega \overline{R}} \setminus D_{\overline{R}}$ with gradient $\leq \varepsilon$. Thus, rescaling back proves the corollary. $\square$

In the next corollary we prove that if the second fundamental form of a compact embedded disk is big at a point, but not necessarily almost its maximum, and the disk satisfies (4.3) in Theorem 4.16 then the disk contains a multi-valued graph, possibly around another point.

**Corollary 4.18.** *Given $N \in \mathbb{Z}_+$, $\omega > 1$, $\varepsilon > 0$, and $T > 0$ there exist $C_1 = C_1(N, \omega, \varepsilon) > 0$, $C_2 = C_2(N, \omega, \varepsilon) > 0$, and $l_1 = l_1(N, \omega, \varepsilon) > 1$ such that the* following holds. If $\Sigma \subset \mathbb{R}^3$ is a compact embedded disk with $0 \in \Sigma \subset B_l(0)$, $\partial\Sigma \subset \partial B_l(0)$, $|A|(0) = C$,

$$|H| \leq \min(C_2, \frac{1}{2l}) \text{ and } [H]_\alpha \leq T$$

for some $l > l_1$ then there exist $p \in \Sigma$, $\overline{R} < \frac{1}{\omega}$, and $0 < \delta < 1$ such that, after a translation that takes $p$ to the origin and possibly after a rotation, $\Sigma$ contains an $N$-valued graph $\Sigma_g$ over $D_{\delta\omega\overline{R}} \setminus D_{\delta\overline{R}}$ with gradient $\leq \varepsilon$.

*Proof.* If

$$\sup_\Sigma |A|^2 \leq 4C^2$$

then we are done, because of Theorem 4.16.

If instead, $\sup_\Sigma |A|^2 > 4C^2$, consider the non-negative function

$$F(x) = (l - |x|)^2 |A(x)|^2.$$

The function $F$ is zero on the boundary of $\Sigma$ and therefore it reaches its maximum at a point in the interior. Let $p$ be such point, i.e.

$$F(p) = \max_\Sigma F(x) = (|p| - l)^2 |A|^2(p) \geq F(0) > 4l^2 C^2.$$

Let $2\sigma < l - |p|$ such that

$$4\sigma^2 |A|^2(p) = 4l^2 C^2.$$

Since $F$ achieves its maximum at $p$,

$$(4.4) \quad \sup_{B_\sigma(p)\cap\Sigma} \sigma^2 |A|^2 \leq \sup_{B_\sigma(p)\cap\Sigma} \sigma^2 \frac{F(x)}{(|x|-l)^2} \leq$$

$$\leq \frac{4\sigma^2}{(|p|-l)^2} \sup_{B_\sigma(p)\cap\Sigma} F(x) = \frac{4\sigma^2}{(|p|-l)^2} F(p) = 4\sigma^2 |A|^2(p).$$

From the weak convex hull property we know that $B_\sigma(p) \cap \Sigma$ consists of a collection of disks. Rescale $B_\sigma(p) \cap \Sigma$ by a factor of $\frac{l}{\sigma} \geq 1$ and translate $p$ to the origin. Let $\Sigma' \subset B_l(0)$ be the rescaled connected component that contains the origin. It follows that $\Sigma'$ is a compact embedded disk such that $\Sigma' \subset \mathcal{B}_l(0)$, $\partial\Sigma' \subset \partial\mathcal{B}_l(0)$,

$$|H'| \leq \frac{\sigma}{l}\min(C_2, \frac{1}{2l}) \leq \min(C_2, \frac{1}{2l}), [H']_\alpha \leq (\frac{\sigma}{l})^{1+\alpha} T \leq T, \sup_{\Sigma'} \leq 4C^2 = 4|A|^2(0).$$

Theorem 4.16 shows that there exists $\overline{R} < \frac{1}{\omega}$ and (after a rotation) an $N$-valued graph $\Sigma_g \subset \Sigma'$ over $D_{\omega\overline{R}} \setminus D_{\overline{R}}$ with gradient $\leq \varepsilon$. Thus, rescaling back proves the corollary. □

## 5. An example illustrating the need for the bound on $\|H\|_\alpha$

In this last section we are going to show that in order to have a multi-valued graph form in a smooth surface, it is not enough to assume that the mean curvature is small relative to the second fundamental form. How small the mean curvature has to be must also depend on the Holder norm of the mean curvature.

Let us consider the graph given by the $C^1$ function

$$u(x,y) = xy \log \sqrt{x^2 + y^2}$$

for $(x, y) \neq (0, 0)$, $u(0, 0) = 0$ and $\nabla u(0, 0) = (0, 0)$.

Let $r = \sqrt{x^2 + y^2}$. One computes,

- $\nabla u = (y \log r + \frac{x^2 y}{r^2}, x \log r + \frac{y^2 x}{r^2})$
- $\Delta u = 4\frac{xy}{r^2}$
- $u_{xy} = \log r + 1 - 2\frac{x^2 y^2}{r^4}$

Recall that the mean curvature $H$ of a graph is given by,

(5.5) $$H = \operatorname{div}\left(\frac{\nabla u}{\sqrt{1 + |\nabla u|^2}}\right).$$

Rewriting equation (5.5) in polar coordinates and computing the mean curvature of $u$ near the origin we obtain

(5.6) $$H(r \cos \theta, r \sin \theta) = 4 \cos \theta \sin \theta + o(r).$$

Let us take a sequence of mollifications $u_\sigma(x, y)$ of $u$ over a disk and, for any $C > 0$ let

$$\overline{u}_\sigma^C = \frac{C}{\max Hess(u_\sigma)} u_\sigma$$

so that

$$\max Hess(\overline{u}_\sigma^C) = C$$

independently on $\sigma$. A computation shows that the sequence of graphs given by the functions $\overline{u}_\sigma^C$ converges $C^1$ to a plane. In fact, because the hessian has been normalized, the graphs in the sequence have second fundamental forms that have size $C$ at the origin. Furthermore, the supremum of the mean curvature is going to zero. Nonetheless, being graphs, these surfaces do not contain a multi-valued graph around any point.

It might appear that this contradicts Corollary 4.18, but this is not the case: Using equations (5.6) one can show that near the origin, as $\sigma$ goes to infinity, the mean curvature of $\overline{u}_\sigma^C$ behaves like the function

$$4C \frac{\cos \theta \sin \theta}{\log r} + o(\frac{r}{\log r}).$$

Since

$$\lim_{r \to 0} r^\alpha \log r = 0 \quad \text{for any } \alpha > 0$$

the $C^\alpha$ norms of the mean curvatures of these graphs are not uniformly bounded with respect to $\sigma$.

## References

[1] G. P. Bessa, L. P. Jorge, and G. Oliveira-Filho. Half-space thorems for minimal surfaces with bounded curvature. *J. Differential Geom.*, 37:493–508, 2001.

[2] T. H. Colding and W. P. Minicozzi II. The Calabi-Yau conjectures for embedded surfaces. To appear in Annals of Math., Preprint math. DG/0404197 (2004).

[3] T. H. Colding and W. P. Minicozzi II. The space of embedded minimal surfaces of fixed genus in a 3-manifold V; Fixed genus. Preprint math.DG/0509647 (2005).

[4] T. H. Colding and W. P. Minicozzi II. Complete properly embedded minimal surfaces in $\mathbb{R}^3$. *Duke Math. J.*, 107:421–426, 2001. MR1823052, Zbl 1010.49025.

[5] T. H. Colding and W. P. Minicozzi II. Embedded minimal disks: proper versus non-proper - global versus local. *Transactions of A.M.S.*, 356(1):283–289, 2003. MR2020033, Zbl 1046.53001.

[6] T. H. Colding and W. P. Minicozzi II. The space of embedded minimal surfaces of fixed genus in a 3-manifold I; Estimates off the axis for disks. *Annals of Math.*, 160:27–68, 2004. MR2119717, Zbl 1070.53031.

[7] T. H. Colding and W. P. Minicozzi II. The space of embedded minimal surfaces of fixed genus in a 3-manifold II; Multi-valued graphs in disks. *Annals of Math.*, 160:69–92, 2004. MR2119718, Zbl 1070.53032.

[8] T. H. Colding and W. P. Minicozzi II. The space of embedded minimal surfaces of fixed genus in a 3-manifold III; Planar domains. *Annals of Math.*, 160:523–572, 2004. MR2123932, Zbl 1076.53068.

[9] T. H. Colding and W. P. Minicozzi II. The space of embedded minimal surfaces of fixed genus in a 3-manifold IV; Locally simply-connected. *Annals of Math.*, 160:573–615, 2004. MR2123933, Zbl 1076.53069.

[10] K. Ecker. A local monotonicity formula for mean curvature flow. *Ann. of Math.*, 154(2):503–525, 2001.

[11] D. Hoffman and W. H. Meeks III. The strong halfspace theorem for minimal surfaces. *Invent. Math.*, 101:373–377, 1990. MR1062966, Zbl 722.53054.

[12] R. López. Area monotonicity for spacelike surfaces with constant mean curvature. *J. Geom. Phys.*, 52(3):353–363, 2003.

[13] W. H. Meeks III, J. Pérez, and A. Ros. Minimal surfaces in $\mathbb{R}^3$ of finite genus. Preprint.

[14] W. H. Meeks III and H. Rosenberg. The uniqueness of the helicoid. *Annals of Math.*, 161:723–754, 2005. MR2153399, Zbl pre02201328.

[15] W. H. Meeks III and H. Rosenberg. The minimal lamination closure theorem. *Duke Math. Journal*, 133(3):467–497, 2006.

[16] W. H. Meeks III and M. Weber. Bending the helicoid. Preprint, arXiv math.DG/0511387.

[17] R. Osserman. *A Survey of Minimal Surfaces*. Dover Publications, New York, 2nd edition, 1986. MR0852409, Zbl 0209.52901.

[18] J. Pérez and A. Ros. Properly embedded minimal surfaces with finite total curvature. In *The Global Theory of Minimal Surfaces in Flat Spaces*, pages 15–66. Lecture Notes in Math 1775, Springer-Verlag, 2002. G. P. Pirola, editor. MR1901613, Zbl 1028.53005.

[19] L. Simon. Lectures on geometric measure theory. In *Proceedings of the Center for Mathematical Analysis*, volume 3, Canberra, Australia, 1983. Australian National University. MR0756417, Zbl 546.49019.

[20] G. Tinaglia. Multi-valued graphs in embedded constant mean curvature disks. *Trans. Amer. Math. Soc.*, 359:143–164, 2007.

[21] F. Xavier. Convex hulls of complete minimal surfaces. *Math. Ann.*, 269:179–182, 1984. MR0759107, Zbl 0528.53009.